\documentclass[12pt,reqno]{amsart}
\usepackage[margin=1.25in]{geometry}        % See geometry.pdf to learn the layout options. There are lots.
\usepackage{graphicx}
\usepackage{amssymb}
\usepackage{aliascnt}
\usepackage{doi}
\usepackage{epstopdf}
\usepackage{esint}

\usepackage[dvipsnames]{xcolor}
\usepackage{tikz,pgfplots}
\usetikzlibrary{patterns.meta}
\usepgflibrary{shadings}
\usepgflibrary {shadings}
\usetikzlibrary {arrows.meta}

\definecolor{ballcolor}{HTML}{53AFDC}
\definecolor{ballcolor1D}{HTML}{84C5E6} % lighter for 1D diagram
\definecolor{twoptcolor}{HTML}{FECF49}
\definecolor{simplexcolor}{HTML}{4CBB51}
\definecolor{mygray}{HTML}{EEEEEE} % lighter gray
\definecolor{mymidgray}{HTML}{A5A5A5} % darker gray

\newtheorem{theorem}{Theorem}[section]

\newaliascnt{corx}{thmx}

\aliascntresetthe{corx}

\newaliascnt{lemma}{theorem}

\aliascntresetthe{lemma}

\newaliascnt{proposition}{theorem}

\aliascntresetthe{proposition}

\newaliascnt{corollary}{theorem}
\newtheorem{corollary}[corollary]{Corollary}
\aliascntresetthe{corollary}

\newaliascnt{conjecture}{theorem}

\aliascntresetthe{conjecture}

\newaliascnt{question}{theorem}

\aliascntresetthe{question}

\theoremstyle{definition}
\newtheorem*{definition*}{Definition}
\newtheorem*{example*}{Example}
\newtheorem*{continuityexamples*}{Continuity examples}
\newtheorem*{remarkmoments}{Remark on minimal moments}
\newtheorem*{remarkvariational}{Remark on variational capacity}

\newcommand{\B}{{\mathbb B}}
\newcommand{\e}{\varepsilon}

\newcommand{\R}{{\mathbb R}}
\newcommand{\Rn}{{{\mathbb R}^n}}

\newcommand{\Sph}{{\mathbb S}}

\DeclareMathOperator{\sinc}{sinc}

\DeclareMathOperator{\dist}{dist}
\DeclareMathOperator{\diam}{diam}

\DeclareMathOperator{\capzero}{Cap_0}
\DeclareMathOperator{\capone}{Cap_1}

\DeclareMathOperator{\capd}{Cap_\mathit{d}}

\DeclareMathOperator{\capp}{Cap_\mathit{p}}
\DeclareMathOperator{\capq}{Cap_\mathit{q}}
\DeclareMathOperator{\capntwo}{Cap_{\mathit{n}-2}}

\newcommand{\arxiv}[1]{%
\href{https://arxiv.org/abs/#1}{ArXiv:#1}}

\title{Riesz capacity: monotonicity, continuity, diameter and volume}
\author{Carrie Clark and Richard S. Laugesen}
\email{carriec2@illinois.edu, Laugesen@illinois.edu}
\address{University of Illinois, Urbana, IL 61801, USA}

\begin{document}

\begin{abstract}
Properties of Riesz capacity are developed with respect to the kernel exponent $p \in (-\infty,n)$, namely that capacity is monotonic as a function of $p$, that its endpoint limits recover the diameter and volume of the set, and that capacity is left-continuous with respect to $p$ and is right-continuous provided (when $p \geq 0$) that an additional hypothesis holds. Left and right continuity properties of the equilibrium measure are obtained too. 
\end{abstract}

\maketitle

\vspace*{-12pt}

\section{\bf Introduction and results}
\label{sec:intro}

\subsection*{Overview}
Riesz capacity measures the size of a set in Euclidean space through an interaction energy with kernel $1/|x-y|^p$. Electrostatic capacity in $3$ dimensions has $p=1$, for example, while the logarithmic kernel in the plane arises as $p \to 0$.

This paper develops properties of Riesz capacity as a function of $-\infty<p<n$, for sets in $n$ dimensional space. We show:
\begin{itemize}
\item[\small $\bullet$] monotonicity of capacity with respect to $p$ (\autoref{th:Rieszmonotonicity}),
\item[\small $\bullet$] diameter and volume as endpoint cases (\autoref{th:Rieszmonotonicity} and \autoref{co:Rieszvolume}),
\item[\small $\bullet$] left-continuity of capacity with respect to $p$, and right-continuity provided an additional hypothesis holds when $p \geq 0$ (\autoref{th:continuity}), 
\item[\small $\bullet$] logarithmic capacity of a nice set equals the limit of its Riesz capacity as $p \to 0$ (\autoref{co:loglimit}), 
\item[\small $\bullet$] analogous left- and right-continuity of equilibrium measures (\autoref{th:eqmmeasure}). 
\end{itemize}

\subsection*{Related work} Inspiration for the volume endpoint case ($p \nearrow n$) comes from work of Hardin and Saff \cite[Theorem 9.3.3]{BHS19} on the discrete charge problem with $p=n$: in their result, the volume provides the dominant term in the energy as the number of point charges tends to infinity. A further motivation is Calef and Hardin's result \cite{CH09} that the Riesz equilibrium measure converges to Lebesgue measure as $p \nearrow n$. Our approach is quite different from either of these motivating works, though, relying instead on a characterization of Riesz capacity due to G\"{o}tz \cite{G03}.

Continuity of capacity and equilibrium measure have antecedents in work of Laugesen \cite{L93,L22} and Pouliasis \cite{P11,P21b} in which the set rather than the kernel is deformed. 

Probabilistic approaches have been employed by Watanabe \cite[p.{\,}489]{W83}, Betsakos \cite{B04a,B04b} and M\'{e}ndez--Hern\'{a}ndez \cite{MH06} when $p \in [n-2,n)$ to minimize Riesz capacity among sets of given volume. Such methods could perhaps also yield monotonicity and continuity results, which in this paper we tackle via purely analytic methods. 

\subsection*{Logarithmic, Riesz, and Newtonian energy and capacity}
Consider a compact set $K$ in $\Rn, n \geq 1$.  The \textbf{logarithmic energy} of $K$ is
\[
V_{log}(K) = \min_\mu \int_K \! \int_K \log \frac{1}{|x-y|} \, d\mu(x) d\mu(y)
\]
where the minimum is taken over all probability measures on $K$, that is, positive unit Borel measures. For the empty set we define $V_{log}(\emptyset)=+\infty$. The energy is greater than $-\infty$ since $|x-y|$ is bounded. If the energy is less than $+\infty$ then a unique measure attains the minimum, called the logarithmic equilibrium measure. (For claims here and below about equilibrium measures, see the references in \autoref{sec:background}.)

The \textbf{Riesz $p$-energy} of $K$ is
\[
V_p(K) = 
\begin{cases}
\min_\mu \int_K \! \int_K |x-y|^{-p} \, d\mu(x) d\mu(y) , & p > 0 , \\
\max_\mu \int_K \! \int_K |x-y|^{-p} \, d\mu(x) d\mu(y) , & p < 0 , 
\end{cases}
\]
where the minimum and maximum are taken over all probability measures on $K$. For the empty set we define $V_p(\emptyset)$ to equal $+\infty$ when $p>0$ and equal $0$ when $p<0$. When $p>0$, the energy is either positive or $+\infty$, and if the energy is finite then the minimum is attained by a unique equilibrium measure. When $p<0$, the energy is automatically finite, and the equilibrium measure (attaining the maximum) is unique if $K$ is nonempty and $-2<p<0$. Notice we do not define Riesz energy when $p=0$. 

The \textbf{Riesz $p$-capacity} is 
\begin{equation} \label{eq:capdef}
\capp(K) = 
\begin{cases}
V_p(K)^{-1/p} , & p \neq 0 , \\
\exp(-V_{log}(K)) , & p = 0 .
\end{cases}
\end{equation}
When $p \geq 0$, the capacity is positive if and only if the energy is finite. When $p<0$, the capacity is positive provided $K$ contains more than one point. The $0$-capacity is also known as \textbf{logarithmic capacity}. The \textbf{Newtonian} energy $V_{n-2}(K)$ and capacity $\capntwo(K)$ are the special case where $n \geq 3$ and $p=n-2$. 

In cases where the Riesz equilibrium measure is unique, we call it the $p$-equilibrium measure of $K$. The $0$-equilibrium measure is simply logarithmic equilibrium measure. 

Capacity is clearly monotonic with respect to set inclusion, with $K_1 \subset K_2$ implying $\capp(K_1) \leq \capp(K_2)$. And capacity scales linearly, with 
\[
\capp(sK)=s \capp(K) , \qquad s>0 .
\]

\subsection*{Justification of the capacity definition.} The definition of Riesz capacity in \eqref{eq:capdef} is the most natural choice, in our opinion, because it makes capacity decreasing with respect to $p$, as \autoref{th:Rieszmonotonicity} will show. 

Other authors, such as Landkof \cite{L72} (for $0<p<n)$ and Borodachov, Hardin and Saff \cite{BHS19} (for $-\infty<p<n)$ define capacity to be $V_p(K)^{-1}$. That definition does not result in monotonicity with respect to $p$. Also, the definition in \cite{BHS19} for $p$-energy when $p<0$ is the negative of ours, which creates an abrupt change in behavior at $p=0$. 

Note that Landkof uses $n-p$ as the Riesz exponent instead of $p$, when $0<p<n$, and so his ``$p$'' must change when the ambient dimension $n$ changes --- even if the set and the energy under consideration remain the same. For example, the $1$-energy of a planar set would arise under Landkof's convention from the choice ``$p$''$=2-1$, but if the planar set is regarded as lying in $3$ dimensions then the same energy arises from taking ``$p$''$=3-1$. The definitions in this paper avoid such changes in exponent: the $p$-capacity of $K$ is unaffected by the dimension $n$ in which the set is regarded as residing. 

\begin{figure}
\begin{center}
\includegraphics[width=0.5\textwidth]{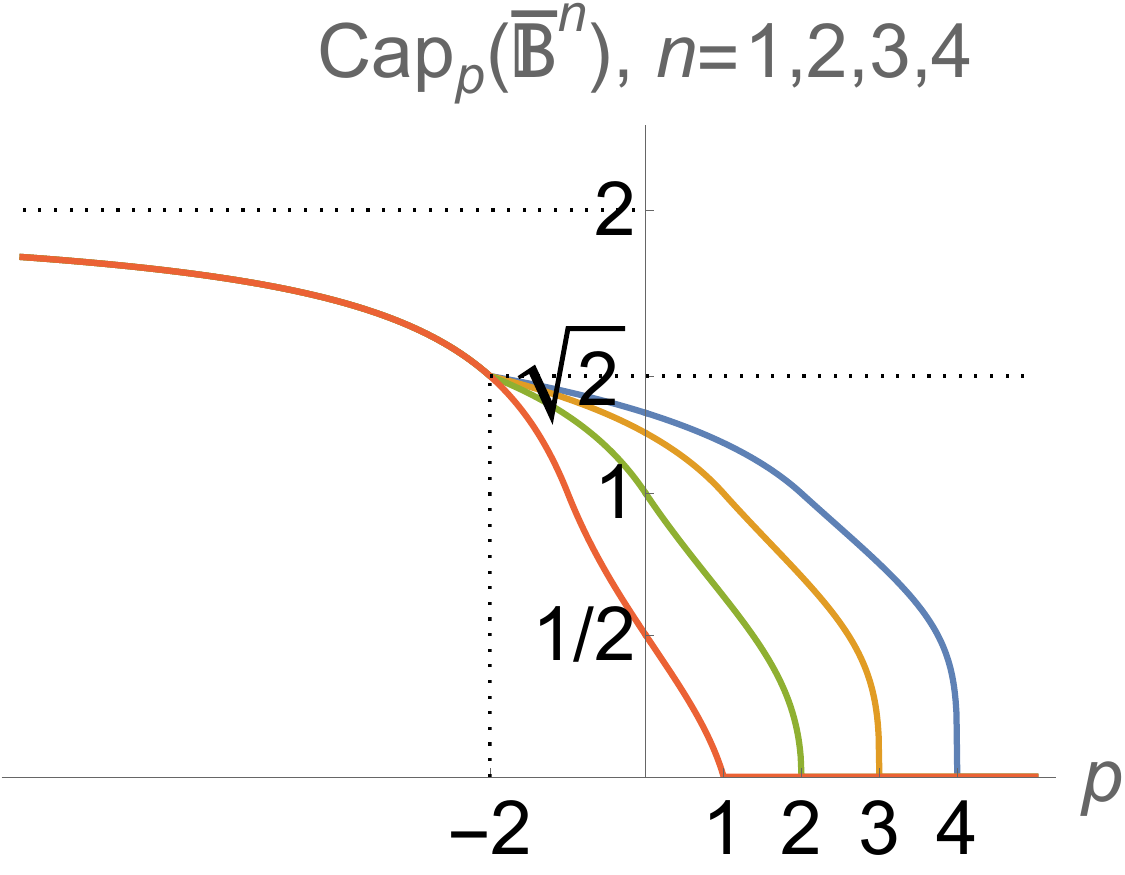}
\end{center}
\caption{\label{fig:ballplots} Monotonicity and continuity of the Riesz capacity $\capp(\overline{\B}^n)$ of the closed unit ball  as a function of $p$, for four values of $n$. The plots rely on formulas from \autoref{sec:background}. The intercept at $p=0$ is the logarithmic capacity $\capzero(\overline{\B}^n)$, with values $1/2, 1, 2/\sqrt{e}, \sqrt[4]{e}$ for $n=1,2,3,4$ respectively. When $p \leq -2$, the capacity of the ball is $2^{1+1/p}$ for each $n$, with limiting value $2$ (the diameter of the ball) as $p \to -\infty$. The horizontal line at height $\sqrt{2}$ is the upper envelope of the family of curves, because for each fixed $p>-2$, one can show using Stirling's formula that $\lim_{n \to \infty} \capp(\overline{\B}^n) = \sqrt{2}$. 
}
\end{figure}

The advantages of definition \eqref{eq:capdef} are illustrated by \autoref{fig:ballplots}, which plots the capacity of the unit ball as a continuous, decreasing function of $p \in \R$. 

\subsection*{Monotonicity and continuity of Riesz capacity wrt $p$}
\label{sec:cont}

Fundamental properties of Riesz capacity as a function of $p$ are stated in the next result and proved in \autoref{sec:proofmonotonicity}. Write $m_n$ for Lebesgue measure on $\Rn$. 
\begin{theorem}[Monotonicity, and limits wrt $p$] \label{th:Rieszmonotonicity}
Let $K \subset \Rn, n \geq 1$, be compact and nonempty.  

(a) (Monotonicity) $\capp(K)$ is a decreasing function of $p \in \R$, and is strictly decreasing on the interval where it is positive. 

(b) (Diameter: $p \to -\infty$) For $p <0$, one has $2^{1/p} \leq \capp(K)/\diam(K) \leq 1$. Hence 
\[
\lim_{p \to -\infty} \capp(K) = \diam(K) .
\]
In fact, when $n=1$ one has $\capp(K) = 2^{1/p} \diam(K)$ for all $p \leq -1$.

(c) (Volume: $p \to n$) 
\[
\capp(K) \sim \left( \frac{m_n(K)}{|\Sph^{n-1}|} (n-p) \right)^{\! \! 1/p} \qquad \text{as $p \nearrow n$.}
\]

(d) $\capp(K)=0$ for all $p \geq n$. 
\end{theorem}
\begin{corollary}[Volume as limiting case of capacity] \label{co:Rieszvolume}
If $K \subset \Rn, n \geq 1$, is compact then 
%If $K$ is a compact subset of $\Rn, n \geq 1$, then 
\[
m_n(K) = |\Sph^{n-1}| \lim_{p \nearrow n} \frac{\capp(K)^p}{n-p} .
\]
\end{corollary}
Parts (a) and (c) of the theorem are new to the best of our knowledge. Part (b) is an easy consequence of the definitions.  

Regarding part (d), much more is known: for $K \subset \Rn$ with Hausdorff dimension $d>0$, one has $\capp(K)>0$ for all $p < d$, while $\capp(K)=0$ for all $p>d$; see \cite[Theorems 4.3.1 and 4.3.3]{BHS19}. In the borderline case $p=d$, if the Hausdorff $d$-measure of $K$ is finite then $\capd(K)=0$. Since the compact set $K$ in \autoref{th:Rieszmonotonicity} has dimension $d \leq n$ and finite $n$-dimensional Lebesgue measure, it follows in particular that $\capp(K)=0$ for $p \geq n$. 

\autoref{co:Rieszvolume} is extended to lower dimensional subsets of $\Rn$, such as $d$-dimensional submanifolds and other ``strongly rectifiable'' sets of integer dimension, by Fan and Laugesen \cite{FL24}. They show that Riesz capacity recovers the $d$-dimensional Hausdorff measure of the set as $p \nearrow d$. 

The next result says that Riesz capacity is always left-continuous with respect to $p$, and is right-continuous too under an additional hypothesis.  
\begin{theorem}[Continuity properties wrt $p$] \label{th:continuity}
Let $K$ be a compact subset of $\Rn, n \geq 1$. 

(i) $\capp(K)$ is continuous for $p<0$.

(ii) $\capp(K)$ is left-continuous for $p \in \R$.

(iii) $\capp(K)$ is right-continuous (and hence continuous) at $p \geq 0$ if $\capp(K)>0$ and for some $p_*>p$ the equilibrium measure $\mu_p$ satisfies 
\begin{equation} \label{eq:pstar}
\int_K \! \int_K \frac{1}{|x-y|^{p_*}} \, d\mu_p d\mu_p < \infty .
\end{equation}
\end{theorem}
This theorem is proved in \autoref{sec:continuity}. The special case at $p=0$ says:
\begin{corollary}[Logarithmic capacity as limit of Riesz capacity] \label{co:loglimit}
If $K$ is a compact subset of $\Rn, n \geq 1$, then
\begin{align*}
\lim_{p \nearrow 0} \capp(K) 
& = \capzero(K) \\
& \geq \lim_{p \searrow 0} \capp(K) .
\end{align*}
Equality holds if $\capzero(K)>0$ and hypothesis \eqref{eq:pstar} is satisfied for some $p_*>0=p$.
\end{corollary}

\begin{continuityexamples*}
1.\ (Balls) The capacity $\capp(B)$ of a closed ball $B$ is continuous at every $p \in \R$, by the explicit formulas in \autoref{sec:background}. See \autoref{fig:ballplots}. 

2.\ (Newtonian) If $K$ is a compact set with smooth boundary in $\Rn, n \geq 2$, then $\capp(K)$ is continuous at the Newtonian exponent $p=n-2$. Indeed, the equilibrium measure $\mu_{n-2}$ on $K$ is known in this case to equal a smooth density function times surface area measure on $\partial K$, from which $\mu_{n-2}$ is easily shown to satisfy hypothesis \eqref{eq:pstar} for $p=n-2<p_*<n-1$. Part (iii) of the theorem now yields continuity at $p=n-2$. 

Thus for smoothly bounded planar sets, \autoref{co:loglimit} implies that logarithmic capacity equals the limit of Riesz capacity as $p \to 0$. 

It seems reasonable to conjecture that $\capp(K)$ is continuous for all $p$, when $K$ has smooth boundary or else satisfies some suitably weaker regularity condition. But we do not know how to verify hypothesis \eqref{eq:pstar} when $p \neq n-2$. The only relevant result we have found in the literature is by Hardin, Reznikov, Saff and Volberg \cite[Corollary 2.9]{HRSV19}, who showed that if $K$ is a smooth $(n-1)$-dimensional submanifold of $\Rn$ and $p \in [n-2,n-1)$ then the equilibrium measure $\mu_p$ is bounded above by a multiple of surface area measure and hence \eqref{eq:pstar} holds for $p_* \in (p,n-1)$. 

3.\ Right continuity can fail. For example, Carleson \cite[Theorem 5 on p.{\,}35]{C67} constructed a compact set in $1$ dimension with positive logarithmic capacity and Hausdorff dimension $0$, so that $\capzero(K)>0$ while $\capp(K)=0$ for all $p>0$. 
%Are there similar counterexamples for each $p>0$? 
\end{continuityexamples*}

\subsection*{Continuous dependence of equilibrium measure wrt $p$}
We finish with some continuous dependence properties of the equilibrium measure, with respect to $p$. The short proof is in \autoref{sec:continuousdependence}. 
\begin{theorem} \label{th:eqmmeasure}
Let $K$ be a compact subset of $\Rn, n \geq 1$, take $q \in \R$, and suppose $\capq(K)>0$. 

(a) (Continuity) If $-2<q<0$ then $\mu_p$ is weak-$*$ continuous at $p=q$. 

(b) (Left-continuity) If $q \geq 0$ then $\mu_p$ is weak-$*$ left-continuous at $p=q$. 

(c) (Right-continuity) If $q \geq 0$ and for some $q_*>q$ the equilibrium measure $\mu_q$ satisfies 
\begin{equation*} 
\int_K \! \int_K \frac{1}{|x-y|^{q_*}} \, d\mu_q d\mu_q < \infty 
\end{equation*}
then $\mu_p$ is weak-$*$ right-continuous at $p=q$. 
\end{theorem}
Calef and Hardin \cite{CH09} showed for strongly rectifiable $d$-dimensional sets that as $p \nearrow d$, the equilibrium measure $\mu_p$ converges to normalized Hausdorff measure. Our left-continuity result in \autoref{th:eqmmeasure}(b) does not apply to their situation when $q=d$, because $\capd(K)=0$ by \cite[Theorem 4.3.1]{BHS19}. 

\begin{remarkmoments} Extremal properties of equilibrium measure are developed in our recent paper \cite{CL24a}, where we show for $p=n-2$ and $p=n-1$ that the moments of $\mu_p$ are minimal when $K$ is a ball, among all sets of given capacity.
\end{remarkmoments}

\begin{remarkvariational} 
Another family of functionals that includes Newtonian capacity is the family of variational capacities, defined by $\min_u \int_\Rn |\nabla u|^q \, dx$ where $u \geq 1$ on $K$ and $u \to 0$ at infinity. The choice $q=2$ gives Newtonian capacity, up to a constant factor. For other values of $q$ there seems to be no direct connection with the Riesz capacities studied in this paper. 
\end{remarkvariational}

\section{\bf Proof of \autoref{th:Rieszmonotonicity} --- monotonicity and limiting values}
\label{sec:proofmonotonicity}

%We may suppose throughout the proof of the theorem that $K$ is nonempty, since the capacity of the empty set is zero for all $p$.

%
\subsection*{Proof of \autoref{th:Rieszmonotonicity}(a): monotonicity}  

Let $\mu$ be a probability measure on $K$. First suppose $0 \leq p < q$. If $p>0$ then Jensen's inequality implies that the $L^p$-norm is less than or equal to the $L^q$-norm, 
\begin{equation} \label{eq:monot1}
\left( \int_K \int_K \frac{1}{|x-y|^p} \, d\mu d\mu \right)^{\! \! 1/p} \leq \left( \int_K \int_K \frac{1}{|x-y|^q} \, d\mu d\mu \right)^{\! \! 1/q} .
\end{equation}
Taking the minimum of the left side over all measures $\mu$ yields that $V_p(K)^{1/p} \leq V_q(K)^{1/q}$, and so $\capp(K) \geq \capq(K)$. To handle $p=0$, observe that  
\begin{equation}  \label{eq:monot2}
\exp \left( q \int_K \int_K \log \frac{1}{|x-y|} \, d\mu d\mu \right) \leq \int_K \int_K \frac{1}{|x-y|^q} \, d\mu d\mu
\end{equation}
by Jensen's inequality. Taking the $q$-th root and minimizing on the left with respect to $\mu$ gives $\exp \left( V_{log}(K) \right) \leq V_q(K)^{1/q}$, and so $\capzero(K) \geq \capq(K)$. 

Thus $\capp(K)$ is a decreasing function of $p \geq 0$. 

Now suppose $p<q \leq 0$, so that $-p > -q \geq 0$. If $-q>0$ then Jensen's inequality implies   
\begin{equation}  \label{eq:monot3}
\left( \int_K \int_K |x-y|^{-p} \, d\mu d\mu \right)^{\! \! -1/p} \geq \left( \int_K \int_K |x-y|^{-q} \, d\mu d\mu \right)^{\! \! -1/q} .
\end{equation}
Taking the maximum over $\mu$ on the left side gives that $\capp(K) \geq \capq(K)$. Lastly, to handle $-q=0$ we note 
\begin{equation}  \label{eq:monot4}
\int_K \int_K |x-y|^{-p} \, d\mu d\mu \geq \exp \left( p \int_K \int_K \log \frac{1}{|x-y|} \, d\mu d\mu \right) 
\end{equation}
by Jensen's inequality. Taking the $(-p)$-th root and maximizing on the left with respect to $\mu$ gives $V_p(K)^{-1/p} \geq \exp \left( -V_{log}(K) \right)$, and so $\capp(K) \geq \capzero(K)$. 

Hence $\capp(K)$ is a decreasing function when $p \leq 0$. 

To prove the strict monotonicity statement in the theorem, suppose $\capp(K)=\capq(K)>0$ in the argument above. We will deduce a contradiction. Consider a $q$-equilibrium measure $\mu$ on $K$. (That measure might not be unique when $q \leq -2$, but uniqueness is not needed for the following argument.) Notice the right sides of the inequalities \eqref{eq:monot1}--\eqref{eq:monot4} are all positive and finite, since the $q$-capacity of $K$ is positive, and so the left sides are all positive and finite too. Equality holds in each inequality, due to our assumption that $\capp(K)=\capq(K)$, and so the equality conditions in Jensen's inequality imply in each case that $|x-y|$ is constant $(\mu \times \mu)$-a.e. 

Write $C$ for that constant value, so that $C = \capq(K)>0$. Let $0<R<C/2$ and cover the compact set $K$ with finitely many balls of radius $R$. At least one of those balls $B$ has positive measure, $\mu(B)>0$. The product $B \times B$ has positive $\mu \times \mu$ measure, but $|x-y| \leq 2R < C$ when $(x,y) \in B \times B$, contradicting the construction of $C$ and thus completing the proof of strict monotonicity. 

\subsection*{Proof of \autoref{th:Rieszmonotonicity}(b): the diameter}
Let $p<0$. If $x,y \in K$ then $|x-y|^{-p} \leq \diam(K)^{-p}$, and so after integrating with respect to $d\mu(x) d\mu(y)$ for an arbitrary probability measure $\mu$ on $K$ we obtain the upper bound $V_p(K) \leq \diam(K)^{-p}$. To get a lower bound, by choosing points $a,b \in K$ that achieve the diameter and letting $\mu = (\delta_a + \delta_b)/2$ be a sum of point charges at those points, we find
\[
V_p(K) \geq \int_K \! \int_K |x-y|^{-p} \, d\mu d\mu = \frac{1}{2} \diam(K)^{-p} .
\]
Combining the two inequalities and taking the $(-p)$-th root implies
\[
\diam(K) \geq \capp(K) \geq 2^{1/p} \diam(K) , \qquad p<0 .
\]
Hence $\capp(K) \to \diam(K)$ as $p \to -\infty$. 

Now suppose $n=1$ and let $a=\min K$ and $b=\max K$, so that $K \subset [a,b]$. Then 
\[
\capp(K) \leq \capp([a,b])=2^{1/p} (b-a) , \qquad p \leq -1 ,
\]
by rescaling the formula $\capp([-1,1])=2^{1+1/p}$ from \autoref{sec:background}. Thus $\capp(K) \leq 2^{1/p} \diam(K)$, which along with the earlier lower bound means that equality holds. 

\subsection*{Proof of \autoref{th:Rieszmonotonicity}(c): showing $\geq$ as $p \nearrow n$}
If $K$ has zero Lebesgue measure, $m_n(K)=0$, then \eqref{eq:Cappnliminf} below holds trivially. So suppose $m_n(K)>0$ and let $\mu=(m_n|_K)/m_n(K)$ be Lebesgue measure restricted to $K$ and then normalized to have unit measure. Using this $\mu$ as a trial measure for the Riesz $p$-energy gives that 
\[
V_p(K) \leq \frac{1}{m_n(K)^2} \int_K \int_K \frac{1}{|x-y|^p} \, dx dy 
\]
when $0<p<n$. Note $K \subset \overline{\B}^n(y,R)$ where $R$ is the diameter of $K$, and so   
\begin{align*}
V_p(K) 
& \leq \frac{1}{m_n(K)^2} \int_K \int_{\overline{\B}^n(R)} \frac{1}{|z|^p} \, dz dy \\
& = \frac{1}{m_n(K)} \frac{ |\Sph^{n-1}| R^{n-p}}{n-p} 
\end{align*}
by spherical coordinates. Hence $\limsup_{p \nearrow n} (n-p) V_p(K) \leq |\Sph^{n-1}| / m_n(K)$, or equivalently 
\begin{equation} \label{eq:Cappnliminf}
\liminf_{p \nearrow n} \frac{\capp(K)}{(n-p)^{1/p}} \geq \left( \frac{m_n(K)}{|\Sph^{n-1}|} \right)^{\! \! 1/n} ,
\end{equation}
which gives the $\geq$ direction of \autoref{th:Rieszmonotonicity}(c). 

\subsubsection*{Alternative proof of inequality \eqref{eq:Cappnliminf}.} Riesz energy can alternatively be characterized when $p>0$ as follows \cite[formula (3) on p.{\,}65]{G03}: 
\begin{equation} \label{eq:Rieszalt}
V_p(K) = p \min_\mu \int_0^\infty \! \int_K \mu(\B^n(x,r)) \, d\mu(x) \, r^{-p-1} \, dr 
\end{equation}
where the minimum is taken over probability measures on $K$. If $V_p(K)$ is finite then the unique minimizer is the $p$-equilibrium measure on $K$.

To apply this characterization, we take $\mu$ as above to be Lebesgue measure restricted to $K$ and normalized to have unit measure. Then \eqref{eq:Rieszalt} yields that 
\begin{align*}
V_p(K) 
& \leq \frac{p}{m_n(K)^2} \int_0^\infty \! \int_K m_n(\B^n(x,r) \cap K) \, dx \, r^{-p-1} \, dr \\
& \leq \frac{p}{m_n(K)^2} \left( \int_0^1 \! \int_K m_n(\B^n(x,r)) \, dx \, r^{-p-1} \, dr + \int_1^\infty \! \int_K m_n(K) \, dx \, r^{-p-1} \, dr \right) \\
& = \frac{|\Sph^{n-1}|}{m_n(K)} \frac{p}{n} \frac{1}{n-p} + 1
\end{align*}
Hence again $\limsup_{p \nearrow n} (n-p) V_p(K) \leq |\Sph^{n-1}| / m_n(K)$, which gives \eqref{eq:Cappnliminf}.

\subsection*{G\"{o}tz's formula for Riesz energy} 
To prove the other direction of \autoref{th:Rieszmonotonicity}(c), we will rely in the proof below on a more sophisticated characterization of Riesz energy that is due to G\"{o}tz \cite[formula (2) on p.{\,}65]{G03}:  
\begin{equation} \label{eq:Gotz}
V_p(K) = A(p,n) \min_\mu \int_0^\infty \! \int_\Rn \mu(\B^n(x,r))^2 \, dx \, r^{-n-p-1} \, dr 
\end{equation}
where $p>0$ and the minimum is taken over probability measures on $K$. 
%If $V_p(K)$ is finite then the minimum is achieved if and only if $\mu$ is the $p$-equilibrium measure on $K$. 
Here the constant is defined to satisfy $1/A(p,1) = 2^{p+1}/p(p+1)$  when $n=1$ and for $n \geq 2$ it satisfies  
\begin{align*}
\frac{1}{A(p,n)} 
& = \frac{2^{p+1}}{p} m_{n-1}(\B^{n-1}) \int_0^1 (1-r^2)^{(n-1)/2} r^p \, dr \\
& = \frac{2^p}{p} \frac{\pi^{(n-1)/2}}{\Gamma(\frac{n+1}{2})} B \! \left(\frac{n+1}{2},\frac{p+1}{2}\right)
%& = \frac{2^p}{p(p+n)} \int_{\Sph^{n-1}} |\xi_1|^p \, d\xi \label{eq:Apn}
\end{align*}
by substituting $r=\sin \theta$ in order to obtain the beta function \cite[{\S}5.12]{DLMF}. By expressing the beta function in terms of gamma functions,  we find when $p=n$ that 
\[
\frac{1}{A(n,n)} = \frac{2^n}{n^2} \frac{\pi^{(n-1)/2}  \Gamma(\frac{n+1}{2})}{\Gamma(n)} = \frac{2\pi^{n/2}}{n^2 \Gamma(\frac{n}{2})} = \frac{|\Sph^{n-1}|}{n^2} ,
\]
where the second equality relies on the duplication formula for the gamma function. Hence 
\begin{equation}
A(n,n) m_n(\B^n)^2 = |\Sph^{n-1}| . \label{eq:Annformula}
\end{equation}

\subsection*{Proof of \autoref{th:Rieszmonotonicity}(c): showing $\leq$ as $p \nearrow n$}
We aim to apply G\"{o}tz's formula to obtain an upper bound on the $\limsup$ of the capacity as $p \nearrow n$. 

Let $\e>0$ and write 
\[
K(\e) = K + \overline{\B^n(\e)}
\]
for the compact set of points within distance $\e$ of $K$. By discarding the parts of G\"{o}tz's energy formula \eqref{eq:Gotz} with either $r \geq \e$ or $x \notin K(\e)$, we get a lower bound  
\begin{align*}
V_p(K) 
& \geq A(p,n) \inf_\mu \int_0^\e \! \int_{K(\e)} \mu(\B^n(x,r))^2 \, dx \, r^{-n-p-1} \, dr \\
& \geq A(p,n) \inf_\mu \int_0^\e \frac{1}{m_n(K(\e))} \left( \int_{K(\e)} \mu(\B^n(x,r)) \, dx \right)^{\! \! 2} r^{-n-p-1} \, dr 
\end{align*}
by Cauchy--Schwarz. When $r < \e$ and $x \in \Rn \setminus K(\e)$, the ball $\B^n(x,r)$ does not intersect $K$ and so $\mu(\B^n(x,r))=0$. Hence the inside integral evaluates to 
\[
\int_\Rn \mu(\B^n(x,r)) \, dx = \mu(K) m_n(\B^n(r)) = m_n(\B^n) r^n ,
\]
from which we deduce the estimate 
\[
V_p(K) \geq \frac{A(p,n) m_n(\B)^2}{m_n(K(\e))} \frac{\e^{n-p}}{n-p} .
\]
Therefore
\begin{align*}
\liminf_{p \nearrow n} \, (n-p) V_p(K) 
& \geq \frac{A(n,n) m_n(\B)^2}{m_n(K(\e))} \\
& \to \frac{|\Sph^{n-1}|}{m_n(K)} \qquad \text{as $\e \to 0$,}
\end{align*}
where \eqref{eq:Annformula} was used in the numerator.  Hence
\begin{equation} \label{eq:Cappnlimsup}
\limsup_{p \nearrow n} \frac{\capp(K)}{(n-p)^{1/p}} \leq  \left( \frac{m_n(K)}{|\Sph^{n-1}|} \right)^{\! \! 1/n} .
\end{equation}
This $\limsup$ estimate in \eqref{eq:Cappnlimsup} combined with the earlier $\liminf$ inequality in \eqref{eq:Cappnliminf} establishes the asymptotic limit for \autoref{th:Rieszmonotonicity}(c). 

\subsection*{Proof of \autoref{co:Rieszvolume}} The corollary is equivalent to \autoref{th:Rieszmonotonicity}(c).

\subsection*{Proof of \autoref{th:Rieszmonotonicity}(d)} \autoref{co:Rieszvolume} shows that $\capp(K) \to 0$ as $p \nearrow n$, and hence $\capp(K) = 0$ for all $p \geq n$ by monotonicity in part (a). 

\subsection*{Remark} Parts (a) and (b) of the theorem could be extended to Riesz-type capacities in a metric space with energy kernel $\dist(x,y)^{-p}$. The proofs of parts (c) and (d) depend on properties of Euclidean space, though, and so any extension of those parts would require additional hypotheses on the metric space.

\section{\bf Proof of \autoref{th:continuity} --- continuity properties}
\label{sec:continuity}

We may assume $K$ is nonempty, since the empty set has $\capp(\emptyset)=0$ for all $p$. 

\subsection*{Proof of \autoref{th:continuity}(i): $p<0$} Continuity of the energy $V_p(K)$ as a function of $p<0$ is clear, since the kernel $|x-y|^{-p}$ is continuous with respect to $x,y,p$ jointly and $K$ is compact. Hence the capacity $\capp(K)$ is continuous when $p<0$. 

Let us explain this argument in more detail, since the technically simple case when $p$ is negative sheds light on the later proof for positive $p$. Fix $q<0$. We will prove that $p \mapsto \capp(K)$ is continuous at $p=q$. Consider an arbitrary sequence $p_j \to q < 0$ with $p_j<0$ for all $j$. For each $j$, choose an equilibrium measure $\mu_j$ on $K$, so that  
\[
\int_K \! \int_K |x-y|^{-p_j} \, d\mu_j(x) d\mu_j(y) = V_{p_j}(K) .
\]
(It does not matter whether this equilibrium measure is unique.) By the Helly selection principle, after passing to a subsequence we may suppose that $\mu_j$ converges weak-$*$ to some probability measure $\mu$ on $K$. 

Weak-$*$ convergence implies weak-$*$ convergence of the product measures (see \cite[Lemma 6.4]{W81}), and so 
\begin{align}
V_q(K) 
& \geq \int_K \! \int_K |x-y|^{-q} \, d\mu d\mu \qquad \text{by definition of the energy for $q<0$} \label{eq:negequality} \\
& = \lim_{j \to \infty} \int_K \! \int_K |x-y|^{-q} \, d\mu_j d\mu_j \qquad \text{by weak-$*$ convergence} \notag \\
& = \lim_{j \to \infty} \int_K \! \int_K |x-y|^{-p_j} \, d\mu_j d\mu_j \qquad \text{since $p_j \to q$} \notag \\
& = \lim_{j \to \infty} V_{p_j}(K) . \notag
\end{align}
The original sequence $p_j$ was arbitrary and so $V_q(K) \geq \limsup_{p \to q} V_p(K)$, which is equivalent to 
\[
\capq(K) \geq \limsup_{p \to q} \capp(K) .
\] 
(Remember that $p$ and $q$ are negative in this part of the argument.)

To get an inequality in the reverse direction, let $\mu_q$ be an equilibrium measure for $V_q(K)$ and observe for $p<0$ that
\begin{align*}
V_p(K) 
& \geq \int_K \! \int_K |x-y|^{-p} \, d\mu_q d\mu_q \qquad \text{by definition of the energy for $p<0$} \\
& \to \int_K \! \int_K |x-y|^{-q} \, d\mu_q d\mu_q \qquad \text{as $p \to q$} \\
& = V_q(K) .
\end{align*}
Hence 
\[
\liminf_{p \to q} \capp(K) \geq \capq(K) .
\] 
Combining the two inequalities shows $\lim_{p \to q} \capp(K) = \capq(K)$, giving continuity at $p=q<0$. Also, the fact that equality holds when the two inequalities are combined means that equality holds in \eqref{eq:negequality}, and so the weak-$*$ limiting measure $\mu$ in the proof above is a $q$-equilibrium measure for $K$. 

\subsection*{Proof of \autoref{th:continuity}(ii): left-continuity} Let $q \geq 0$. We will prove that $p \mapsto \capp(K)$ is left continuous at $p=q$. Consider an arbitrary sequence $p_j \to q$, which for now we do not restrict to approaching from the left. For each $j$, choose an equilibrium measure $\mu_j$ on $K$, so that  
\[
\int_K \! \int_K |x-y|^{-p_j} \, d\mu_j(x) d\mu_j(y) = V_{p_j}(K) .
\]
(If the energy equals $+\infty$ then any measure $\mu_j$ will do.) After passing to a subsequence we may suppose that $\mu_j$ converges weak-$*$ to some probability measure $\mu$ on $K$. 

Suppose first that $q>0$, and fix $N>0$. Weak-$*$ convergence implies weak-$*$ convergence of the product measures, and so 
\begin{align*}
& \int_K \! \int_K \min \left( N , |x-y|^{-q} \right) d\mu d\mu \\
& = \lim_{j \to \infty} \int_K \! \int_K \min \left( N , |x-y|^{-q} \right) d\mu_j d\mu_j \\
& = \lim_{j \to \infty} \int_K \! \int_K \min \left( N , |x-y|^{-p_j} \right) d\mu_j d\mu_j \qquad \text{since $p_j \to q$} \\
& \leq \liminf_{j \to \infty} \int_K \! \int_K |x-y|^{-p_j} \, d\mu_j d\mu_j  \qquad \text{by dropping the ``$\min$''} \\
& = \liminf_{j \to \infty} V_{p_j}(K) 
\end{align*}
by our choice of $\mu_j$ as an equilibrium measure for $p=p_j$. Letting $N \to \infty$ on the left side, we conclude from monotone convergence and the definition of the energy as a minimum that 
\begin{equation} \label{eq:qliminf}
V_q(K) \leq \int_K \! \int_K |x-y|^{-q} \, d\mu d\mu \leq \liminf_{j \to \infty} V_{p_j}(K) .
\end{equation}
The sequence $p_j$ was arbitrary, and so $V_q(K) \leq \liminf_{p \to q} V_p(K)$, which is equivalent to 
\[
\limsup_{p \to q} \capp(K) \leq \capq(K) .
\] 

When $p$ increases to $q$ we further have $\lim_{p \nearrow q} \capp(K) \geq \capq(K)$ by the decreasing property of capacity in \autoref{th:Rieszmonotonicity}, and so 
\[
\lim_{p \nearrow q} \capp(K) = \capq(K) ,
\]
which is left-continuity. 

For later use, note that left-continuity forces equality to hold throughout \eqref{eq:qliminf} when $p_j \nearrow q$, which implies in that case that $\mu$ is a $q$-equilibrium measure for $K$. 

Suppose next that $p_j \to q=0$, and fix $N>0$. Require $p_j \neq 0$ for all $j$. By the weak-$*$ convergence, 
\begin{align*}
& \int_K \! \int_K \log \min \left( N , \frac{1}{|x-y|} \right) d\mu d\mu \\
& = \lim_{j \to \infty} \int_K \! \int_K \log \min \left( N , \frac{1}{|x-y|} \right) d\mu_j d\mu_j \\
& = \lim_{j \to \infty} \frac{1}{p_j} \log \left[ 1 + p_j \int_K \! \int_K \log \min \left( N , \frac{1}{|x-y|} \right) d\mu_j d\mu_j \right] \\
& = \lim_{j \to \infty} \frac{1}{p_j} \log \left[ \int_K \! \int_K \exp \left( p_j \log \min \left( N , \frac{1}{|x-y|} \right) \right) d\mu_j d\mu_j \right] 
\end{align*}
by the exponential series, noting here that $\min ( N , 1/|x-y| )$ is bounded both above and below and so the second order term in the series is bounded by $O(p_j)^2$, independently of $x, y \in K$. The last limit increases when we drop the $\min$ (regardless of whether $p_j$ is positive or negative), and so 
\begin{align*}
& \int_K \! \int_K \log \min \left( N , \frac{1}{|x-y|} \right) d\mu d\mu \\
& \leq \liminf_{j \to \infty} \frac{1}{p_j} \log \int_K \! \int_K \exp \left( p_j \log \frac{1}{|x-y|} \right) d\mu_j d\mu_j \\
& = \liminf_{j \to \infty} \frac{1}{p_j} \log \int_K \! \int_K \frac{1}{|x-y|^{p_j}} \, d\mu_j d\mu_j \\
& = \liminf_{j \to \infty} \frac{1}{p_j} \log V_{p_j}(K)
\end{align*}
by choice of $\mu_j$ as a $p_j$-equilibrium measure. Letting $N \to \infty$ on the left side, applying monotone convergence, and recalling the definition of the logarithmic energy as a minimum, we find 
\begin{equation} \label{eq:zeroliminf}
V_{log}(K) \leq \int_K \! \int_K \log \frac{1}{|x-y|} \, d\mu d\mu \leq \liminf_{j \to \infty} \log V_{p_j}(K)^{1/p_j} .
\end{equation}
The sequence $p_j$ was arbitrary, and so $V_{log}(K) \leq \liminf_{p \to 0} \log V_p(K)^{1/p}$, which by exponentiating is equivalent to $\limsup_{p \to 0} \capp(K) \leq \capzero(K)$. When $p$ increases to $0$ we also know $\lim_{p \nearrow 0} \capp(K) \geq \capzero(K)$ by the decreasing property of capacity, and so $\lim_{p \nearrow 0} \capp(K) = \capzero(K)$, which gives left-continuity at $q=0$. 

Notice that left-continuity forces equality to hold throughout \eqref{eq:zeroliminf} when $p_j \nearrow 0$, and so in that case $\mu$ is a logarithmic equilibrium measure for $K$.

\subsection*{Proof of \autoref{th:continuity}(iii): right-continuity at $q \geq 0$} Continue with the notation from part (ii) above. Let $q \geq 0$. The task is to prove right-continuity of $\capp(K)$ at $p=q$. 

Suppose $\capq(K)>0$ and that for some $q_*>q$ the equilibrium measure $\mu_q$ satisfies 
\begin{equation} \label{eq:qbound}
\int_K \! \int_K \frac{1}{|x-y|^{q_*}} \, d\mu_q d\mu_q < \infty .
\end{equation}
Hypothesis \eqref{eq:qbound} guarantees that $V_{q_*}(K)$ is finite and so $\operatorname{Cap_\mathit{q_*}}(K)>0$. Capacity decreases with the Riesz exponent, and so for all $p \leq q_*$ we have $\capp(K)>0$. Thus the equilibrium measure $\mu_p$ exists and is unique for all $0 \leq p \leq q_*$. 

Let $p_j \searrow q$ be an arbitrary sequence with $p_j<q_*$ for all $j$, and denote by $\mu_j$ the equilibrium measure corresponding to $p_j$. The Helly selection principle ensures after passing to a subsequence that $\mu_j$ converges weak-$*$ to a probability measure $\mu$ on $K$. 

Suppose first that $q>0$. (Later we treat $q=0$.) Inequality \eqref{eq:qliminf} says that 
\begin{equation} \label{eq:lastminute1}
V_q(K) \leq \int_K \! \int_K \frac{1}{|x-y|^q} \, d\mu d\mu \leq \liminf_{j \to \infty} V_{p_j}(K) .
\end{equation}
Using $\mu_q$ as a trial measure for $V_{p_j}(K)$ yields that 
\begin{align*}
\limsup_{j \to \infty} V_{p_j}(K) 
& \leq \limsup_{j \to \infty} \int_K \! \int_K \frac{1}{|x-y|^{p_j}} \, d\mu_q d\mu_q \\
& = \int_K \! \int_K \frac{1}{|x-y|^q} \, d\mu_q d\mu_q = V_q(K)
\end{align*}
by dominated convergence, with the dominator provided by \eqref{eq:qbound} since $p_j<q_*$. Combining the last two inequalities shows that $\lim_{j \to \infty} V_{p_j}(K) = V_q(K)$. The original sequence $p_j$ was arbitrary and so 
\[
\lim_{p \searrow q} V_p(K) = V_q(K) ,
\]
which gives right-continuity of capacity, namely $\lim_{p \searrow q} \capp(K) = \capq(K)$. 

The argument shows that equality holds in \eqref{eq:lastminute1} and so $\int_K \! \int_K 1/|x-y|^q \, d\mu d\mu = V_q(K)$. Hence $\mu$ must be the equilibrium measure $\mu_q$. 

\subsection*{Proof of \autoref{th:continuity}(iii): right-continuity at $q=0$} 
Now take $q=0$, so that $q_*> p_j \searrow 0$. Using the logarithmic equilibrium measure $\mu_0$ as a trial measure for $V_{p_j}(K)$ yields that 
\begin{equation} \label{eq:lastminute2}
\log V_{p_j}(K)^{\! 1/p_j} 
\leq \frac{1}{p_j} \log \int_K \! \int_K \frac{1}{|x-y|^{p_j}} \, d\mu_0 d\mu_0 .
\end{equation}
Let $x,y \in K, x \neq y$. Taylor's theorem for the function $p \mapsto1/ |x-y|^p$ says that
\[
\frac{1}{|x-y|^p} - 1 - p \log \frac{1}{|x-y|} = \frac{p^2}{2} \left( \log \frac{1}{|x-y|} \right)^{\! \! 2} \frac{1}{|x-y|^{p(x,y)}} 
\]
for some number $p(x,y) \in (0,p)$. Suppose $0<p<q_*/2$, so that $p(x,y)<q_*/2$. Then the right side of the last formula is bounded by $Ap^2/|x-y|^{q_*}$ for some constant $A$, and that quantity is integrable with respect to $d\mu_0 d\mu_0$ by  hypothesis \eqref{eq:qbound}. Thus
\[
\int_K \! \int_K \frac{1}{|x-y|^{p_j}} \, d\mu_0 d\mu_0 = 1 + p_j V_{log}(K) + p_j^2 Q_j
\]
where $Q_j$ is bounded independently of $j$, and so from \eqref{eq:lastminute2} we deduce
\begin{align*}
\limsup_{j \to \infty} \log V_{p_j}(K)^{\! 1/p_j} 
& \leq \limsup_{j \to \infty} \frac{1}{p_j} \log \left( 1 + p_j V_{log}(K) + p_j^2 Q_j \right) \\
& = V_{log}(K) .
\end{align*}

Combining this inequality with the upper bound \eqref{eq:zeroliminf} on the logarithmic energy shows that $\lim_{j \to \infty} \log V_{p_j}(K)^{\! 1/p_j} = V_{log}(K)$. Since the original $p_j$ sequence was arbitrary, we conclude that 
\[
\lim_{p \searrow 0} \log V_p(K)^{\! 1/p} = V_{log}(K) .
\]
Exponentiating yields $\lim_{p \searrow 0} \capp(K) = \capzero(K)$, which is right-continuity at $q=0$. 

Lastly, this argument using \eqref{eq:zeroliminf} shows $\int_K \int_K \log 1/|x-y| \, d\mu d\mu = V_{log}(K)$, so that $\mu$ must equal the logarithmic equilibrium measure $\mu_0$.

\section{\bf Proof of \autoref{th:eqmmeasure}}
\label{sec:continuousdependence}

The theorem assumes that the capacity is positive at $p=q$, and hence for all smaller $p$-values too. Thus the equilibrium measure is unique for each $p$ and $q$ in the following proof, noting that the hypotheses of part (a) restrict to $q> -2$ in order to ensure uniqueness. 

In the proof of \autoref{th:continuity}(i) we showed that an arbitrary sequence converging to $q$ has a subsequence $p_j$ such that $\mu_{p_j} \to \mu_q$ as $j \to \infty$. It follows that 
\[
\mu_p \to \mu_q \text{\ \ weak-$*$ \quad as $p \to q < 0$,}
\]
which gives the continuity in \autoref{th:eqmmeasure}(a). The proof of \autoref{th:continuity}(ii) implies similarly that 
\[
\mu_p \to \mu_q \text{\ \ weak-$*$ \quad as $p \nearrow q \geq 0$,}
\]
giving left-continuity in part (b) of the theorem. Under the hypotheses of part (c), the proof of \autoref{th:continuity}(iii) shows likewise that  
\[
\mu_p \to \mu_q \text{\ \ weak-$*$ \quad as $p \searrow q \geq 0$,}
\]
which is the desired right-continuity in the theorem.

\section*{Acknowledgments}
Richard Laugesen's research was supported by grants from the Simons Foundation (\#964018) and the National Science Foundation ({\#}2246537).

\appendix

\section{\bf Potential theoretic background}
\label{sec:background}

\subsection*{Existence and uniqueness of equilibrium measure} 

A weak-$*$ compactness argument (Helly selection principle) applied to probability measures on $K$ shows the existence of a measure achieving the minimum in the definition of the energy, for both the logarithmic ($p=0$) and Riesz cases ($p \neq 0$); see \cite[Lemma 4.1.3]{BHS19}. Existence of an equilibrium measure is easy when $p<0$, of course, since the Riesz kernel is continuous in that case. 

Uniqueness of the logarithmic equilibrium measure ($p=0$) in all dimensions is treated by Borodachov, Hardin and Saff \cite[Theorem 4.4.8]{BHS19}. That theorem also provides uniqueness of the equilibrium measure for $-2<p<0$. When $p>0$, uniqueness can be found in \cite[Theorem 4.4.5]{BHS19}. Remember when $p \geq 0$ that to get uniqueness of the equilibrium measure one assumes $K$ has finite energy, in other words, positive capacity. 

%Uniqueness when $0<p<n$ is of course a standard result, for example in Landkof \cite[pp.{\,}132--133]{L72}. Uniqueness of the logarithmic equilibrium measure in $2$ dimensions is also well known: for example, \cite[pp.\,133,167--168]{L72}, or the appealing variant in Saff and Totik \cite[Theorem I.1.3, Lemma I.1.8]{ST97}. 

\subsection*{Riesz capacities and equilibrium measures for the unit ball} In what follows, remember from \autoref{sec:intro} that the Riesz $p$-capacity in this paper differs by a $p$-th root from the definition employed in \cite{BHS19}, and that we use Riesz exponent $p$ rather than $n-p$ as in \cite{L72}. Under our definition, the $p$-capacity of the interval ($n=1$) is  
\begin{equation*} 
\capp(\overline{\B}^1) = 
\begin{cases}
2^{1+1/p} & \text{for $p \leq -1$,} \\
1/2 & \text{for $p=0$,} \\
\left( \frac{\Gamma(1/2)}{\Gamma((1-p)/2)\Gamma(1+p/2)} \right)^{\! \! 1/p} & \text{for $-1<p<1, \ p \neq 0$,} 
\end{cases}
\end{equation*}
by \cite[Proposition 4.6.1]{BHS19}, with unique equilibrium measure given by $\mu=(\delta_{-1}+\delta_{+1})/2$ when $p \leq -1$; and when $-1<p<1$ the equilibrium measure is 
\[
d\mu(x) = \frac{1}{B(1/2,(1+p)/2)} \, \frac{dx}{(1-x^2)^{(1-p)/2}} , \qquad -1<x<1 ,
\]
with normalizing constant 
\[
B(1/2,(p+1)/2)= \frac{\Gamma(1/2) \Gamma((p+1)/2)}{\Gamma(1+p/2)} .
\]

Assume now that $n \geq 2$. The Riesz $p$-capacity of the unit ball is  
\begin{equation} \label{eq:Rieszball}
\capp(\overline{\B}^n) = 
\begin{cases}
2^{1+1/p} & \text{for $p \leq -2$,} \\
2 \exp \! \left( \frac{\psi((n-1)/2) - \psi(n-1)}{2} \right) & \text{for $p=0$,} \\
2 \left( \frac{\Gamma(n-1-p/2) \Gamma((n-1)/2)}{\Gamma((n-1-p)/2) \Gamma(n-1)} \right)^{\! \! 1/p}
& \text{if $-2<p \leq n-2, \ p \neq 0$,} \\
\left( \frac{\Gamma(n/2)}{\Gamma((n-p)/2)\Gamma(1+p/2)} \right)^{\! \! 1/p} & \text{if $n-2<p<n$,} 
%\left( \frac{\Gamma(1/2)}{\Gamma((1-p)/2)\Gamma(1+p/2)} \right)^{\! \! 1/p} & \text{for $0<p<1=n$,}
\end{cases}
\end{equation}
by \cite[Theorem 4.5.8, Proposition 4.6.1, Proposition 4.6.4, Theorem 4.6.6, Theorem 4.6.7]{BHS19}; alternatively, for some of the cases one may consult Landkof \cite[p.{\,}163]{L72}. Here $\psi=(\log \Gamma)^\prime$ is the digamma function. Some formulas have been rewritten using standard identities for the gamma function. 
%(In Landkof's formulas, replace his dimension $p$ with $n$, and replace his exponent $\alpha$ with $n-p$. Also, when $0<p \leq n-2$ we manipulated Landkof's expression using the duplication formula \cite[5.5.5]{DLMF} with $z=n/2$ and the functional equation $\Gamma(z+1)=z\Gamma(z)$, in order to arrive at the expression above.) 

Continue assuming $n \geq 2$. The equilibrium measure for $\capp(\overline{\B}^n)$ is unique when $-2<p<n$, being given by normalized $(n-1)$-dimensional surface area measure on $\Sph^{n-1}$ when $-2<p \leq n-2$; and when $n-2<p<n$ it is  
\[
d\mu(x) = \frac{2}{|\Sph^{n-1}|} \frac{1}{B \big( \frac{n}{2},\frac{p-(n-2)}{2} \big)} \, \frac{dx}{(1-|x|^2)^{(n-p)/2}} , \qquad |x|<1 ,
\]
with the normalizing constant given by the beta function 
\[
B \left( \frac{n}{2},\frac{p-(n-2)}{2} \right) = \frac{\Gamma(\frac{n}{2}) \Gamma(\frac{p-(n-2)}{2})}{\Gamma(1+\frac{p}{2})} .
\]
When $p < -2$ the equilibrium measure of the ball is not unique: it has the form $\mu=(\delta_{-\xi}+\delta_{+\xi})/2$ for an arbitrary point $\xi \in \Sph^{n-1}$. When $p=-2$, equilibrium measures are those probability measures supported on the unit sphere that have center of mass at the origin. For these claims, see \cite[Theorems 4.6.6, 4.6.7]{BHS19}.

\subsection*{Appealing special cases}
By the formula above, the unit disk ($n=2$) has 
\[
\capp(\overline{\B}^2) = \left( \frac{\sin (p\pi/2)}{p\pi/2} \right)^{\! \! 1/p} = \sinc (p\pi/2)^{1/p} , \qquad p \in (0,2) .
\]
Taking $p=1$, the unit disk has Newtonian capacity $\capone(\overline{\B}^2) = 2/\pi$.  

The unit ball in $3$ dimensions has  
\[
\capp(\overline{\B}^3) = 2(1-p/2)^{1/p}, \qquad p \in (-2,1] ,
\]
by simplifying the formula above with $\Gamma(1+p/2)=(p/2) \Gamma(p/2)$ and then using the reflection formula for the gamma function \cite[(5.5.3)]{DLMF}. 

In $n$-dimensions, the unit ball has Newtonian capacity 
\[
\capntwo(\overline{\B}^n) = 1 , \qquad n \geq 1 , 
\]
since the third formula in \eqref{eq:Rieszball} simplifies to $1$ when $p=n-2, n \geq 3$, with the help of the duplication formula \cite[(5.5.5)]{DLMF}. 

The logarithmic capacity ($p=0$) of the ball in the first few dimensions is 
\[
\capzero(\overline{\B}^n) = 
\begin{cases}
1/2 & \text{for $n=1$,} \\
1 & \text{for $n=2$,} \\
2/\sqrt{e} & \text{for $n=3$,} \\
\sqrt[4]{e} & \text{for $n=4$.}
\end{cases}
\]

\bibliographystyle{plain}

\end{document}